\documentstyle[12pt]{article}
\makeatletter
\renewcommand\thesection{\@Roman\c@section}
\renewcommand\thesubsection{\thesection.\@arabic\c@subsection}
\makeatother

\begin{document}
\begin{titlepage}
\begin{flushright}
math.QA/0001154
\end{flushright}
\vskip.3in

\begin{center}
{\Large \bf Level-one Highest Weight Representation of  
$U_q[\widehat{sl(N|1)}]$ and 
Bosonization of the Multi-component Super $t-J$
Model }
\vskip.3in
{\large Wen-Li Yang $^{a,b}$~~and Yao-Zhong Zhang $^b$}
\vskip.2in
{\em $~^a$ Institute of Modern Physics, Northwest University, 
Xian 710069 ,China \\
$~^b$ Department of Mathematics, University of Queensland, Brisbane,
     Qld 4072, Australia}
\end{center}

\vskip 2cm
\begin{center}
{\bf Abstract}
\end{center}

We study the level-one irreducible highest weight representations of 
the quantum affine superalgebra
$U_q[\widehat{sl(N|1)}]$, and calculate their  characters 
and supercharacters. 
We obtain bosonized q-vertex operators acting on the irreducible 
$U_q[\widehat{sl(N|1)}]$-modules and  derive the exchange relations 
satisfied by the  vertex operators. 
We give the bosonization of the  multi-component super  
$t-J$ model by using the bosonized  vertex operators.

\vskip 3cm
\noindent{\bf Mathematics Subject Classifications (1991):} 17B37, 81R10, 
                 81R50, 16W30

\end{titlepage}


\def\a{\alpha}
\def\b{\beta}
\def\d{\delta}
\def\e{\epsilon}
\def\ve{\varepsilon}
\def\g{\gamma}
\def\k{\kappa}
\def\l{\lambda}
\def\o{\omega}
\def\t{\theta}
\def\s{\sigma}
\def\u{\nu}
\def\D{\Delta}
\def\L{\Lambda}

\def\R{\overline{R}}
\def\S{{sl(N|1)}}
\def\hS{{\widehat{sl(N|1)}}}
\def\hG{{\widehat{gl(N|N)}}}
\def\R{{\cal R}}
\def\hR{{\hat{\cal R}}}
\def\C{{\bf C}}
\def\P{{\bf P}}
\def\Z2{{{\bf Z}_2}}
\def\Z{{\bf Z}}
\def\T{{\cal T}}
\def\H{{\cal H}}
\def\F{{\cal F}}
\def\V{\overline{V}}
\def\trho{{\tilde{\rho}}}
\def\tphi{{\tilde{\phi}}}
\def\tT{{\tilde{\cal T}}}
\def\uS{{U_q[sl(N|1)]}}
\def\uqsnh{{U_q[\widehat{sl(N|1)}]}}
\def\uqgnh{{U_q[\widehat{gl(N|N)}]}}
\def\uq1h{{U_q[\widehat{gl(1|1)}]}}
\def\uqs2h{{U_q[\widehat{sl(2|1)}]}}


\def\beq{\begin{equation}}
\def\eeq{\end{equation}}
\def\bea{\begin{eqnarray}}
\def\eea{\end{eqnarray}}
\def\ba{\begin{array}}
\def\ea{\end{array}}
\def\no{\nonumber}
\def\lt{\left}
\def\rt{\right}
\newcommand{\bq}{\begin{quote}}
\newcommand{\eq}{\end{quote}}

\newtheorem{Theorem}{Theorem}
\newtheorem{Definition}{Definition}
\newtheorem{Proposition}{Proposition}
\newtheorem{Lemma}{Lemma}
\newtheorem{Corollary}{Corollary}
\newcommand{\proof}[1]{{\bf Proof. }
        #1\begin{flushright}$\Box$\end{flushright}}

\newcommand{\sect}[1]{\setcounter{equation}{0}\section{#1}}
\renewcommand{\theequation}{\thesection.\arabic{equation}}

\sect{Introduction}
The purpose of this paper is two-fold. One is to study irreducible 
highest weight representations and q-vertex operators \cite{Fre92} of the 
quantum affine superalgebra $\uqsnh$, $N>2$. Another one is to apply these 
results to bosonize the multi-component super $t-J$ model on an infinite
lattice.

We shall adapt the bosonization technique initiated in \cite{Fre88,Ber89}, 
which turns out to be very powerful in constructing highest weight 
representations and q-vertex operators. Recently, free bosonic
realizations of the level-one representations and ``elementary" q-vertex 
operators have been obtained for $U_q[\widehat{sl(M|N)}]$, $M\neq N$ 
\cite{Kim97} and  $U_q[\widehat{gl(N|N)}]$ \cite{Zha98}. However, these 
free boson representations are not irreducible in general. Moreover, the 
elementary q-vertex operators obtained in \cite{Kim97,Zha98} 
were determined solely from their
commutation relations with the bosonized Drinfeld generators \cite{Dri88} 
of the relevant algebras, and thus one can ask on which representations
these bosonized q-vertex operators act. To construct irreducible 
highest weight representations and q-vertex operators acting on them, 
we need to study in details the structure of the bosonic Fock space
generated by the free boson fields. This has been done for 
$U_q[\widehat{sl(2|1)}]$ \cite{Kim97,Yan98} and  $U_q[\widehat{gl(N|N)}]$, 
$N\leq 2$ \cite{Yan99}. In this paper we treat the 
$U_q[\widehat{sl(N|1)}]$ ($N>2$) case.

Irreducible highest weight representations and bosonized q-vertex
operators acting on them play an essential role in the  
algebraic analysis method of lattice integrable models, which was 
invented by the Kyoto group and  collaborators \cite{Dav93,Jim94}. 
In this approach, the following assumption is the vital key :
\bea
``{ \rm \bf the~ physical~space~ of ~states~of ~the~
model}"=\bigoplus_{\a,\a'}
V(\l_{\a})\otimes V(\l_{\a'})^{*S}\label{Ass}
\eea
where $V(\l_{\a})$ is the level-one irreducible highest weight module 
of the underlying quantum affine algebras and $V(\l_{\a})^{*S}$ is the
dual module 
of $V(\l_{\a})$. By this method, various integrable models have been 
analysed such as   
the higher spin XXZ chains \cite{Idz94,Bou94,Hon98}, the higher rank
cases \cite{Koy94,Dav99}, the twisted $A^{(2)}_2$ case \cite{Hou99}, and 
the face type statistical  models \cite{Luk96,AJMP}.

Spin chain models with quantum superalgebra symmetries have
been the focus of recent studies in the
context of strongly correlated fermion systems 
\cite{Foe93,Ess92,Bra95,Ram96,Pfa96}. It is natural to generalize 
the algebraic analysis method to treat super spin chains on an infinite 
lattice. In \cite{Yan98}, the q-deformed supersymmetric  $t-J$ 
model which has $U'_q[\widehat{sl(2|1)}]$ as its non-abelian symmetry 
has been analysed. However, the super case is
fundamentally different from the non-super case. Unlike the latter,
$U_q[\widehat{sl(2|1)}]$ has infinite number of level-one irreducible 
highest weight representations and the bosonized q-vertex operators act in
all of them. This leads to \cite{Yan98} the assumption that for
the $q$-deformed supersymmetric $t$-$J$ model
$\a$, $\a'$ in (\ref{Ass}) take infinite number of  integer values. 

In this paper we extend the work \cite{Yan98} to treat the 
multi-component $t-J$ model with $U'_q[\widehat{sl(N|1)}]$ ($N>2$) 
symmetry. As we shall see, the level-one irreducible highest weight 
representations of $U_q[\widehat{sl(N|1)}]$ ($N>2$) have similar 
structures as the $N=2$ case. So we shall make the assumption  that 
the physical space of states of the multi-component $t-J$ model 
on an infinite lattice is of the form (\ref{Ass}) with 
$\a,~ \a'$ being any integers.

This paper is organized as follows. After presenting some neccessary
preliminaries, we in section 3 construct the level-one 
irreducible highest weight representations of $\uqsnh$ and calculate 
their (super)characters  by means of the BRST resolution. 
In section 4, we compute the exchange relations of the 
q-vertex operators and show that they form  the 
graded Faddeev-Zamolodchikov algebra. In section 5, we consider the
application of these results to the multi-component  super 
$t-J$ model on an infinite lattice. Generalizing the Kyoto group's work 
\cite{Dav93}, we give  the bosonization of this model  using the
bosonized vertex operators of $\uqsnh$. Finally, 
we compute the one-point correlation functions of the local operators and 
give an integral expression of the correlation functions.

\sect{Preliminaries}

\subsection{ Quantum affine superalgebra $\uqsnh$ }

Let us introduce orthonormal basis $\{\e'_i|i=1,2,\cdots,N+1\}$ with  
the bilinear form $(\e'_i,\e'_j)=\u_i\d_{ij}$, where $\u_i=1$ for
$i\neq N+1$ and $\u_{N+1}=-1$. The classical     fundamental weights are
defined by $\bar{\L}_i=\sum_{j=1}^{i}\e_j$ ($i=1,2,\cdots,N$), 
with  $\e_i=\e'_i -\frac{\u_i}{N-1}\sum_{j=1}^{N+1}\e'_j$.  
Introduce the affine weight $\L_0$ and the null root $\d$ having
$(\L_0,\e'_i)=(\d,\e'_i)=0$ for $i=1,2,\cdots, N+1$ and
$(\L_0,\L_0)=(\d,\d)=0$, $(\L_0,\d)=1$. The affine simple roots
and fundamental weights are given by 
\bea
&& \a_i=\u_i\e'_i-\u_{i+1}\e'_{i+1}, ~~ i=1,2,\cdots, N,~~~~~~
\a_0=\d-\sum_{i=1}^{N}\a_i~,\no\\
&&\L_0=\L_0,~~~~\L_i=\L_0+\bar{\L}_i,~~i=1,2,\cdots,N.
\eea
The Cartan matrix of the affine superalgebra $\hS$ reads as
\bea
(a_{ij})=\left(
\begin{array}{ccccccc}
0 & -1 &  &  & &1\\
-1&  2 &-1&  & & \\
 &-1& 2 & \ddots & &  \\
 && \ddots & \ddots & \ddots &    \\
 &&        &  -1    &       2&-1 \\
1&&        &        &      -1& 0\\
\end{array}
\right) ~~~(i,j=0,1,2,\cdots,N).\label{Cartan-matrix}
\eea

The Quantum affine superalgebra
$\uqsnh$ is a $q$-analogue of the universal
enveloping algebra of $\widehat{sl(N|1)}$ generated by the Chevalley
generators $\{e_i,f_i,q^{h_i},d | i=0,1,2,\cdots,N\}$, where $d$ is
the
usual derivation operator. The $Z_2$-grading
of the generators are
$[e_0]=[f_0]=[e_{N}]=[f_{N}]=1$ and zero otherwise.
The defining relations are
\begin{eqnarray*}
& & [h_i,h_j] =0,~~~ h_id=dh_i,~~~[d,e_i]=\delta_{i,0}e_i,~~~
[d,f_i]=-\delta_{i,0}f_i,\\
& &q^{h_i}e_jq^{-h_i}=q^{a_{ij}}e_j,~~~~
q^{h_i}f_jq^{-h_i}=q^{-a_{ij}}f_j,~~~~
[e_i,f_j] =\delta_{ij}{ q^{h_i}-q^{-h_i} \over q-q^{-1}},\\
&&[e_i,e_j]=[f_i,f_j]=0,~~{\rm for}~~ a_{ij}=0,\\
&&[e_j,[e_j,e_i]_{q^{-1}}]_q=0,~~~[f_j,[f_j,f_i]_{q^{-1}}]_q=0,
~{\rm for}~~ |a_{ij}|=1,~~ j\neq 0,N.
\end{eqnarray*} 
\noindent Here and throughout, $[a,b]_x\equiv ab-(-1)^{[a][b]}x ba$ 
and $[a,b]\equiv [a,b]_1$. We do not write down the  extra  $q$-Serre
relations which can be 
obtained by using Yamane's Dynkin diagram procedure \cite{Yam96}.

$\uqsnh$ is a $\Z_2$-graded quasi-triangular Hopf algebra endowed with
the following coproduct $\D$, counit $\e$ and antipode $S$:
\bea
\D(h_i)&=&h_i\otimes 1+1\otimes h_i,~~~~\D(d)=d\otimes 1+1\otimes d,\no\\
\D(e_i)&=&e_i\otimes 1+q^{h_i}\otimes e_i,~~~~
\D(f_i)=f_i\otimes q^{-h_i}+1\otimes f_i,\no\\
\e(e_i)&=&\e(f_i)=\e(h)=0,\no\\
S(e_i)&=&-q^{-h_i}e_i,~~~~S(f_i)=-f_iq^{h_i},~~~~S(h)=-h,\label{e-s} 
\eea
where $i=0,1,\cdots,N$.
Notice that the antipode $S$ is a $\Z_2$-graded algebra anti-homomorphism.
Namely, for any homogeneous elements $a,b\in\uqsnh$
$S(ab)=(-1)^{[a][b]}S(b)S(a)$, which extends to inhomogeneous elements
through linearity. Moreover,
\beq
S^2(a)=q^{-2\rho}\,a\,q^{2\rho},~~~~\forall a\in\uqsnh,
\eeq
where $\rho$ is an element in the Cartan subalgebra 
such that $(\rho,\a_i)=(\a_i,\a_i)/2$ for any simple root
$\a_i,~i=0,1,2,\cdots,N$. Explicitly,
\beq
\rho=(N-1)d+\bar{\rho}=(N-1)d+\frac{1}{2}\sum_{k=1}^{N}(N-2k)\e'_k
-\frac{1}{2}N\e'_{N+1},\label{P1}
\eeq
which $\bar{\rho}$ is  the half-sum of positive roots
of $sl(N|1)$ .
The multiplication rule on the tensor products is $\Z_2$-graded:
$(a\otimes b)(a'\otimes b')=(-1)^{[b][a']}(aa'\otimes bb')$ for any
homogeneous elements $a,b,a',b'\in \uqsnh$.

$\uqsnh$ can also be realized in terms of the Drinfeld generators
\cite{Dri88} 
$\{X^{\pm,i}_m,~ H^i_n,~ q^{\pm H^i_0}$, $~c,~ d |
m\in {\bf Z},~ n\in{\bf Z}-\{0\},~i=1,2,\cdots,N\}$.  
The $\Z_2$-grading of the Drinfeld generators is
given by $[X^{\pm,N}_m]=1$ for  $m\in{\bf Z}$ and 
zero otherwise.
The relations 
satisfied by the Drinfeld generators read  \cite{Yam96,Zha97}
\bea
&&[c,a]=[d,H^i_0]=[H^i_0, H^{j}_n]=0,~~~~[d,H^i_n]=nH^i_n,
~~~\forall a\in\uqsnh\no\\
&&[d,X^{\pm,i}_n]=nX^{\pm,i}_n,~~~~
q^{H^j_0}X^{\pm, i}_nq^{-H^j_0}=q^{\pm a_{ij}}X^{\pm, i}_n,\no\\
&&[H^i_n, H^{j}_m]=\d_{n+m,0}\frac{[a_{ij}n]_q[nc]_q}{n},~~~~[H^i_n,
X^{\pm, j}_m]=\pm\frac{[a_{ij}n]_q}{n}X^{\pm,j}_{n+m}q^{\pm|n|c/2},\no\\
&&[X^{+,i}_n, X^{-,j}_m]=\frac{\d_{ij}}{q-q^{-1}}\lt(q^{\frac{c}{2}(n-m)}
  \psi^{+,i}_{n+m}-q^{-\frac{c}{2}(n-m)}\psi^{-,i}_{n+m}\rt),\no\\
&&[X^{\pm,i}_n, X^{\pm, j}_m]=0,~~~~{\rm for}~a_{ij}=0,\no\\
&&[X^{\pm,i}_{n+1}, X^{\pm,j}_m]_{q^{\pm a_{ij}}}
  -[X^{\pm,j}_{m+1}, X^{\pm,i}_n]_{q^{\pm a_{ij}}}=0,~~{\rm
for}~~a_{ij}\neq 0,\no\\
&&Sym_{l,m}[X^{\pm,i}_l,[X^{\pm,i}_m,X^{\pm,j}_n]_{q^{-1}}]_q=0,~~
{\rm for }~~a_{ij}=0, i\neq N,\label{drinfeld}
\eea
where 
$\sum_{n\in{\bf Z}}\psi^{\pm,j}_{n}z^{- n}=q^{\pm H^j_0}\exp\lt(
  \pm(q-q^{-1})\sum_{n>0}H^j_{\pm n}z^{\mp n}\rt)$ , and the symbol
$Sym_{k,l}$ means symmetrization with respect to 
$k$ and $l$. We used the standard notation 
 $[x]_q=(q^x-q^{-x})/(q-q^{-1})$. The Chevalley generators are related 
to the Drinfeld generators by the
formulas:
\bea
&&h_i=H_0^i,~~~e_i=X^{+,i}_0,~~~f_i=X^{-,i}_0,~~~
  i=1,2,\cdots,N,~~h_0=c-\sum_{k=1}^{N}H^k_0,\no\\
&&e_0=-[X^{-,N}_0,[X^{-,N-1}_0,\cdots,[X^{-,2}_0,
  X^{-,1}_1]_{q^{-1}}]_{q^{-1}}\cdots]_{q^{-1}}
  \,q^{-\sum_{k=1}^{N}H^k_0},\no\\
&&f_0=q^{\sum_{k=1}^{N}H^k_0}[[\cdots [[X^{+,1}_{-1},
  X^{+,2}_0]_{q}
  ,\cdots,X^{+,N-1}_0]_{q},X^{+,N}_0]_{q}.
\eea

\subsection{Free Bosonic realization of the quantum affine superalgebra 
$\uqsnh$ at level one}

Introduce bosonic oscillators $\{a^i_n,~b_n,~c_n$, $Q_{a^i},Q_{b},
Q_{c}$ $|n\in{\bf Z}, i=1,2,\cdots, N,\}$ which
satisfy the commutation relations
\bea
&&[a^i_n, a^{j}_m]=\d_{n+m,0}\d_{ij}\frac{[n]_q[n]_q}{n},~~~~~
  [a^i_0, Q_{a^{j}}]=\d_{ij},\no\\
&&[b_n, b_m]=-\d_{n+m,0}\frac{[n]_q^2}{n},~~~~~
  [b_0, Q_{b}]=-1,\no\\
&&[c_n, c_m]=\d_{n+m,0}\frac{[n]_q^2}{n},~~~~~
  [c_0, Q_{c}]=1.\label{oscilators}
\eea
The remaining commutation relations are zero. Define
$\{h^i_m|i=1,2,\cdots,N,~m\in\Z\}$ :
\bea
&&h^i_m=a^i_mq^{-|m|/2}-a^{i+1}q^{|m|/2},~~
Q_{h_i}=Q_{a^i}-Q_{a^{i+1}},~~~i=1,2,\cdots,N-1,\no\\
&&h^N_m=a^N_mq^{-|m|/2}+b_mq^{-|m|/2},~~~~
Q_{h_N}=Q_{a^N}+Q_b.
\eea
Let us introduce the notation $h^j(z;\k)=Q_{h_j}+h^j_0\ln z   
   -\sum_{n\neq 0}\frac{h^j_n}{[n]_q}q^{\k |n|}z^{-n}$.
The bosonic fields 
$c(z;\b)$, $b(z;\b)$ and  $h^{*}_j(z;\b)$ are  defined in the same
way. 
Define the Drinfeld currents,
$X^{\pm,i}(z)=\sum_{n\in{\bf Z}}X^{\pm,i}_n z^{-n-1}$, 
$i=1,2,\cdots,N$, and the q-differential operator  
$\partial_zf(z)=\frac{f(qz)-f(q^{-1}z)}{(q-q^{-1})z}$.
Then, the Drinfeld generators of $\uqsnh$  at level one can be 
realized by the free boson fields as \cite{Kim97}
\bea
&&c=1,~~~~H^i_m=h^i_m,~~~~
X^{+,N}(z)=:e^{h^N(z;-\frac{1}{2})}e^{c(z;0)}:
e^{-\sqrt{-1}\pi\sum_{i=1}^{N-1}a^i_0},\no\\
&&X^{-,N}(z)=:e^{-h^N(z;\frac{1}{2})}\partial_z\{e^{-c(z;0)}\}:
e^{\sqrt{-1}\pi\sum_{i=1}^{N-1}a^i_0},\no\\
&&X^{\pm,i}(z)=\pm:e^{\pm
h^i(z;\mp\frac{1}{2})}:e^{\pm\sqrt{-1}\pi a^i_0},~~~~~ i=1,2,\cdots,N-1.
\label{A2}
\eea

\subsection{Bosonization of level-one vertex operators}

In order to construct the vertex operators of $\uqsnh$, we firstly
consider the level-zero representations (i.e. the  evaluation
representations) of $\uqsnh$.

Let $E_{i,j}$ be the $(N+1)\times (N+1)$ matrix whose $(i,j)$-element
is unity and zero elsewhere. Let
$\{v_1,v_2,\cdots,v_{N+1}\}$ be the basis vectors of the (N+1)-dimensional
graded vector space $V$. The $\Z_2$-grading of these basis vectors is
chosen to be $[v_i]=(\u_i+1)/2$.
The $(N+1)$-dimensional level-zero
representation $V_z$ of $\uqsnh$ 
is given  by
\bea
&&e_i=E_{i,i+1},~~~~~~
f_i=\nu_i E_{i+1,i},~~~~~~
t_i=q^{\nu_i E_{i,i}-\nu_{i+1}E_{i+1,i+1}},\no \\
&&e_0=-z E_{N+1,1}, ~~~~
f_0=z^{-1} E_{1,N+1},~~~~
t_0=q^{-E_{1,1}-E_{N+1,N+1}},
\eea
where $i=1,\cdots,N$.
Let $V^{*S}_z$ be the left dual module of $V_z$, 
defined by 
\beq
\pi_{V^{*S}_z}(a)=\pi_{V_z}(S(a))^{st},~~~~\forall a\in\uqsnh,
\eeq
where $st$ denotes the supertansposition.

Now,  we study the level-one
vertex operators \cite{Fre92} of $\uqsnh$. Let
$V(\l)$ be the highest weight $\uqsnh$-module with the highest weight
$\l$ and the highest weight vector $|\l>$. Consider the following
intertwiners of $\uqsnh$-modules
\cite{Jim94}:
\bea
&&\Phi^{\mu V}_\l(z):~~ V(\l)\longrightarrow V(\mu)\otimes V_z,
    ~~~~
\Phi^{\mu V^*}_\l(z):~~ V(\l)\longrightarrow V(\mu)\otimes V^{*S}_z,
     \no\\
&&\Psi^{V\mu}_\l(z):~~ V(\l)\longrightarrow V_z\otimes V(\mu),
    ~~~~
\Psi^{V^*\mu}_\l(z):~~ V(\l)\longrightarrow V^{*S}_z\otimes V(\mu).
     \label{TypeII}
\eea
They are intertwiners in the sense that for any $x\in \uqsnh$
\beq
\Xi(z)\cdot x=\D(x)\cdot\Xi(z),~~~~\Xi(z)=\Phi^{\mu V}_\l(z),~
   \Phi^{\mu V^*}_\l(z),~\Psi^{V\mu}_\l(z),~\Psi^{V^*\mu}_\l(z).
   \label{intertwiner1}
\eeq
We expand the vertex operators as \cite{Jim94}
\bea
&&\Phi^{\mu V}_\l(z)=\sum_{j=1}^{N}\,\Phi^{\mu V}_{\l,j}(z)   
   \otimes v_j,~~~~ 
\Phi^{\mu V^*}_\l(z)=\sum_{j=1}^{N}\,\Phi^{\mu V^*}_{\l,j}(z)
   \otimes v^*_j,\no\\
&&\Psi^{V\mu}_\l(z)=\sum_{j=1}^{N}\,v_j\otimes \Psi^{V\mu}_{\l,j}(z),~~~~
\Psi^{V^*\mu}_\l(z)=\sum_{j=1}^{N}\,v^*_j\otimes \Psi^{V^*\mu}_{\l,j}(z).
\eea
The intertwiners are even,  which implies $[\Phi^{\mu
V}_{\l,j}(z)] = [\Phi^{\mu V^*}_{\l,j}(z)] = [\Psi^{V\mu}_{\l,j}(z)] 
= [\Psi^{V^*\mu}_{\l,j}(z)] = [v_j] = \frac{\u_j+1}{2}$.
According to \cite{Jim94},
$\Phi^{\mu V}_\l(z)~\lt(\Phi^{\mu V^*}_\l(z)\rt)$ is called type I
(dual) vertex operator and
$\Psi^{V\mu}_\l(z)~\lt(\Psi^{V^*\mu}_\l(z)\rt)$ type II (dual) vertex
operator.

Introduce the bosonic operators $\phi_j(z)$, $\phi^*_j(z)$, $\psi_j(z)$
 and $\psi^*_j(z)$ \cite{Kim97}:
\bea
&&\phi_{N+1}(z)=:e^{-h^{*}_N(q^Nz;\frac{1}{2})}
e^{c(q^Nz;0)}(q^Nz)^{\frac{N-2}{2(N-1)}}
:e^{\sqrt{-1}\pi\sum_{i=1}^{N}\frac{1-i}{N-1}a^i_0},
\no\\
&&\nu_l\phi_l(z)(-1)^{[f_l]([v_l]+[v_{l+1}])}
=[\phi_{l+1}(z)~,~f_l]_{q^{\nu_{l+1}}},\no\\
&&\phi^*_1(z)=:e^{h^{*}_{1}(qz;\frac{1}{2})}(q^Nz)^{\frac{N-2}{2(N-1)}} :
e^{-\sqrt{-1}\pi\sum_{i=1}^{N}\frac{1-i}{N-1}a^i_0}
,\no\\
&&-\nu_lq^{\nu_l}\phi^*_{l+1}(z)(-1)^{[f_l]([v_l]+[v_{l+1}])}
=[\phi^*_l(z)~,f_l]_{q^{\nu_l}},\no\\
&&\psi_1(z)=:e^{-h^{*}_{1}(qz;-\frac{1}{2})}(q^Nz)^{\frac{N-2}{2(N-1)}} :
e^{\sqrt{-1}\pi\sum_{i=1}^{N}\frac{1-i}{N-1}a^i_0},\no\\
&&\psi_{l+1}(z)=[\psi_l(z)~,e_l]_{q^{\nu_l}},\no\\
&&\psi^*_{N+1}(z)=:e^{h^{*}_N(q^{2-N}z;-\frac{1}{2})}
\partial_z\{e^{-c(q^{2-N}z;0)}\}(q^Nz)^{\frac{N-2}{2(N-1)}} :
e^{-\sqrt{-1}\pi\sum_{i=1}^{N}\frac{1-i}{N-1}a^i_0},\no\\            
&&-\nu_l\nu_{l+1}q^{-\nu_l}\psi^*_l(z)=[\psi^*_{l+1}(z)~,~e_l]_{q^{\nu_{l+1}}},
\label{Vertex-operator2}
\eea
where 
\begin{eqnarray*}
&&h^{*i}_n=\sum_{j=1}^{N}\frac{[\a_{ij}m]_q[\b_{ij}m]_q}{[(N-1)m]_q[m]_q}h^j_n
~,~~Q^*_{h^i}=\sum_{j=1}^{N}\frac{\a_{ij}\b_{ij}}{N-1}Q_{h^j}
~,~~h^{*i}_0=\sum_{j=1}^{N}\frac{\a_{ij}\b_{ij}}{N-1}h^j,
\end{eqnarray*}
with $\a_{ij}=min(i,j)$, and $\b_{ij}=N-1-max(i,j)$.
Define  the even  operators $\phi(z),\;\phi^*(z),\;\psi(z)$ and
$\psi^*(z)$ by
$\phi(z)=\sum_{j=1}^{N+1}\,\phi_j(z)
   \otimes v_j$, $\phi^*(z)=\sum_{j=1}^{N+1}\,\phi^*_j(z)
   \otimes v^*_j$, $\psi(z)=\sum_{j=1}^{N+1}\,v_j\otimes \psi_j(z)$ and
$\psi^*(z)=\sum_{j=1}^{N+1}\,v^*_j\otimes \psi^*_j(z)$. 
Then the vertex operators $\Phi^{\mu V}_\l(z),\;
\Phi^{\mu V^*}_\l(z),\; \Psi^{V\mu}_\l(z)$ and $\Psi^{V^*\mu}_\l(z)$, if
they exist,
are  bosonized by  $\phi(z),\;\phi^*(z),\;\psi(z)$
and $\psi^*(z)$, respectively \cite{Kim97}.
We remark that our vertex operators differ from those of Kimura et al
\cite{Kim97} by  a  scalar factor $(q^Nz)^{\frac{N-2}{2(N-1)}}$, which  
is needed in order for the vertex operators also 
satisfy (\ref{intertwiner1}) for the element $x=d$.
$\phi(z),\;\phi^*(z),\;\psi(z)$ and $\psi^*(z)$ are referred to as the 
``elementary q-vertex operators" of $\uqsnh$.

\sect{Highest weight $\uqsnh$-modules}   
We begin by defining the Fock module. Denote by
$F_{\l_1,\l_2,\cdots,\l_{N+1};\l_{N+2}}$ the bosonic Fock
space generated by $a_{-m}^i,b_{-m},c_{-m} (m>0)$
 over the vector $|\l_1,\l_2,\cdots,\l_{N+1};\l_{N+2}>$:
\begin{eqnarray*}
F_{\l_1,\l_2,\cdots,\l_{N+1};\l_{N+2}}=
{\bf
C}[a^i_{-1},a^i_{-2},\cdots;b_{-1},b_{-2},\cdots;c_{-1},c_{-2},\cdots]
|\l_1,\l_2,\cdots,\l_{N+1};\l_{N+2}>,
\end{eqnarray*}
where
\begin{eqnarray*}
|\l_1,\l_2,\cdots,\l_{N+1};\l_{N+2}>=
e^{\sum_{i=1}^{N}\l_iQ_{a^i}+\l_{N+1}Q_{b}+\l_{N+2}Q_{c}}|0>.
\end{eqnarray*}
\noindent The vacuum vector $|0>$ is defined by
$a^i_m|0>=b_m|0>=c_m|0>=0$ for $i=1,2,\cdots,N$, and $ m\geq 0$.
Obviously,
\begin{eqnarray*}
&& a^i_m|\l_1,\l_2,\cdots,\l_{N+1};\l_{N+2}>=0,{\rm
~~for}~i=1,2,\cdots,N~{\rm
and}~~m>0,\\
&&b_m|\l_1,\l_2,\cdots,\l_{N+1};\l_{N+2}>= 
c_m|\l_1,\l_2,\cdots,\l_{N+1};\l_{N+2}>=0,~~{\rm for}~
m>0.
\end{eqnarray*}
To obtain the highest weight vectors of $\uqsnh$, we impose the
conditions:
\bea
&&e_i|\l_1,\cdots,\l_{N+1};\l_{N+2}>=0,\;i=0,1,2,\cdots,N,\no\\
&&h_i|\l_1,\cdots,\l_{N+1};\l_{N+2}>=\l^i|\l_1,\cdots,\l_{N+1};\l_{N+2}>
,~~i=0,1,2,\cdots,N.
\eea
Solving these equations, we obtain two classess of solutions:
\begin{enumerate}

\item   
$(\l_1,\cdots,\l_i,\l_{i+1},\cdots,\l_{N+1};\l_{N+2})=(\b+1,\cdots,
\underbrace{\b+1,\;\b}_{i,\;i+1},
\cdots,\b;0)$, where $i=1,\cdots,N$, and $\b$ is arbitrary. It follows
that 
$(\l^0,\l^1,\cdots,\l^i,\l^{i+1}\cdots, \l^{N})$=$
(0,0,\cdots,\underbrace{0,\;1,\;0}_{i-1,\;i,\;i+1},\cdots,0)$ and  we have
the identification
$|\L_i>$=$|\b+1,\cdots,\underbrace{\b+1,\;\b}_{i,\;i+1},\cdots,\b;0>$.

\item   
$(\l_1,\cdots,\l_{N},\l_{N+1};\l_{N+2})=(\b,\cdots,\b,\b-\a;-\a)$, where
$\a$, $\b$ are  arbitrary.  
We have 
$(\l^0,\l^1,\cdots,\l^{N-1}, \l^{N})$=$
(1-\a,0,\cdots,0,\a)$ and 
$|(1-\a)\L_0+\a\L_N>$=$|\b,\cdots,\b,\b-\a;-\a>$.
\end{enumerate}
Associated to the above two classes of solutions are the following
Fock spaces:
\bea
&&\F^{m}_{\b}=\bigoplus_{\{i_1,\cdots,i_N\}\in
\Z}F_{\b+1+i_1,\b+1-i_1+i_2,\cdots,\b+1-i_{m-1}+i_{m},\b-i_{m}+i_{m+1}
,\cdots,\b-i_{N-1}+i_{N},\b+i_N;i_N},\no\\
&&\F_{(\a;\b)}=\bigoplus_{\{i_1,\cdots,i_N\}\in
\Z}F_{\b+i_1,\b-i_1+i_2,\cdots,
\b-i_{N-1}+i_{N},\b-\a+i_N;-\a+i_N},\no
\eea
where $m=1,2,\cdots,N$, and it should be understood that $i_0\equiv 0$. 
However, it is easily seen that 
$\F^{m}_{\b}=F_{(m;\b)}$, $m=1,\cdots,N$. Thus, it is sufficient to
study the Fock space
$\F_{(\a;\b)}$. In the following we shall also restrict ourselves to the
 $\a\in \Z$ case.

It can be shown that the bosonized action of $\uqsnh$ (\ref{A2}) 
on $\F_{(\a;\b)}$
 is closed :
\begin{eqnarray*}
&&\uqsnh \F_{(\a;\b)}=\F_{(\a;\b)}.
\end{eqnarray*}
\noindent Hence each Fack space $\F_{(\a;\b)}$ 
constitutes a  $\uqsnh$-module. However, these modules are  not
irreducible in general. To obtain irreducible subspaces, we
introduce a pair of ghost fields \cite{Kim97}

\begin{eqnarray*}
\eta(z)=\sum_{n\in \Z}\eta_nz^{-n-1}=:e^{c(z)}:,~~~~~~
\xi(z)=\sum_{n\in \Z}\xi_nz^{-n}=:e^{-c(z)}:.
\end{eqnarray*}
\noindent The mode expansion of $\eta(z)$ and $\xi(z)$ is well
defined on  $\F_{(\a;\b)}$ for 
$\a\in \Z$, and  
the modes satisfy the relations 
\bea
&&\xi_m\xi_n+\xi_n\xi_m=\eta_m\eta_n+\eta_n\eta_m=0,~~~~~~
\xi_m\eta_n+\eta_n\xi_m=\delta_{m+n,0}.
\eea
Since $\eta_0\xi_0$ and $\xi_0\eta_0$ qualify as 
projectors,  we use them to decompose $\F_{(\a;\b)}$ into 
a direct sum 
$\F_{(\a;\b)}=\eta_0\xi_0\F_{(\a;\b)}\oplus 
\xi_0\eta_0\F_{(\a;\b)}$  for $\a\in \Z$.  
$\eta_0\xi_0\F_{(\a;\b)}$ is referred to as $Ker_{\eta_0}$ and 
$\xi_0\eta_0\F_{(\a;\b)}=\F_{(\a;\b)}/\eta_0\xi_0\F_{(\a;\b)}$ 
as $Coker_{\eta_0}$. Since $\eta_0$ commutes (or 
anticommutes) with the bosonized action of $\uqsnh$, $Ker_{\eta_0}$ 
and $Coker_{\eta_0}$ are both $\uqsnh$-modules for $\a\in\Z$.
\subsection{Character and supercharacter}

We want to determine  the character and supercharacter formulae of
the $\uqsnh$-modules constructed in the bosonic Fock
space.
We first of all bosonize the derivation operator $d$ as
\bea
d=-\sum_{m\geq 1}\frac{m^2}{[m]_q^2}\{
\sum_{i=1}^{N}h^i_{-m}h^{*i}_{m}
+c_{-m}c_m\}
-\frac{1}{2}\{\sum_{i=1}^{N}h^i_0h^{*i}_0+c_0(c_0+1)\}.
\label{d}
\eea
It  obeys the commutation 
relations
\begin{eqnarray*}
[d,h_i]=0,~~~~[d,h^i_m]=mh^i_m,~~~~[d,X^{\pm,i}_m]=mX^{\pm,i}_m,~~~~
i=1,2,\cdots,N,
\end{eqnarray*}
\noindent as required. Moreover,  $[d,\xi_0]=[d,\eta_0]=0$. 

The character and supercharacter of a $\uqsnh$-module $M$ are defined by  
\bea
Ch_{M}(q;x_1,x_2,\cdots,x_N)&=&tr_M(q^{-d}x_1^{h_1}x_2^{h_2}\cdots 
 x_N^{h_N}
),\no\\
Sch_{M}(q;x_1,x_2,\cdots,x_N)&=&Str_M(q^{-d}x_1^{h_1}x_2^{h_2}\cdots
 x_N^{h_N})\no\\
&=&tr_M((-1)^{N_f}q^{-d}x_1^{h_1}x_2^{h_2}\cdots
 x_N^{h_N}),
\eea
respectively. The Fermi-number operaor $N_f$ can be bosonized as 
\bea 
N_f=\left\{
\begin{array}{ll}
(N-1)b_0&~~~if~~N~even,~i.e.~~N=2L\\
L(\sum_{k=1}^{N}a^i_0-b_0)+c_0&~~~if~~N~odd,~i.e.~~N=2L+1
\end{array}
\right..
\eea
\noindent Indeed, $N_f$ satisfies 
\begin{eqnarray*}
(-1)^{N_f}\Theta(z)=(-1)^{[\Theta(z)]}\Theta(z)(-1)^{N_f},
\end{eqnarray*}
\noindent where $\Theta(z)=X^{\pm,i}(z)$, $\phi_i(z)$, $\phi^*_i(z)$, 
$\psi_i(z)$ and $\psi^*_i(z)$.

We calculate the characters and
supercharacters  by using the BRST resolution
\cite{Yan98}. Let us define the Fock spaces, for $l\in \Z$
\begin{eqnarray*}
\F^{(l)}_{(\a;\b)}=\bigoplus_{\{i_1,\cdots,i_N\}\in
\Z}F_{\b+i_1,\b-i_1+i_2,\cdots,
,\b-i_{N-1}+i_{N},\b-\a+i_N;-\a+i_N+l}.
\end{eqnarray*}
We have $\F^{(0)}_{(\a;\b)}=\F_{(\a;\b)}$. It can be shown that
$\eta_0$ and $\xi_0$ intertwine these Fock spaces as follows:
\begin{eqnarray*}
&&\eta_0:~\F^{(l)}_{(\a;\b)}\longrightarrow
\F^{(l+1)}_{(\a;\b)},~~~~
\xi_0:~\F^{(l)}_{(\a;\b)}\longrightarrow
\F^{(l-1)}_{(\a;\b)}.
\end{eqnarray*}
We have the following  BRST complexes:
\bea
\begin{array}{ccccccc}
\cdots &\stackrel{Q_{l-1}=\eta_0}{\longrightarrow}&
\F^{(l)}_{(\a;\b)}&\stackrel{Q_{l}=\eta_0}{\longrightarrow}&
\F^{(l+1)}_{(\a;\b)}&\stackrel{Q_{l+1}=\eta_0}{\longrightarrow}&
\cdots\\
&&|{\bf O}&&|{\bf O}&&\\
\cdots &\stackrel{Q_{l-1}=\eta_0}{\longrightarrow}&
\F^{(l)}_{(\a;\b)}&\stackrel{Q_{l}=\eta_0}{\longrightarrow}&
\F^{(l+1)}_{(\a;\b)}&\stackrel{Q_{l+1}=\eta_0}{\longrightarrow}&
\cdots
\end{array}
\label{brst}
\eea
where ${\bf O}$ is an operator such that
$\F^{(l)}_{(\a;\b)}\longrightarrow \F^{(l)}_{(\a;\b)}$.
Noting the fact that
$\eta_0\xi_0+\xi_0\eta_0=1$, and
$\eta_0\xi_0$ ($\xi_0\eta_0$) is the projection operator
from $\F^{(l)}_{(\a;\b)}$ to $Ker_{Q_{l}}$
($Coker_{Q_{l}}$),
we  get
\bea
&&Ker_{Q_l}=Im_{Q_{l-1}}, ~~{\rm for~ any~} l\in \Z,
\no\\
&&tr({\bf O})|_{Ker_{Q_l}}=tr({\bf O})|_{Im_{Q_{l-1}}}
=tr({\bf O})|_{Coker_{Q_{l-1}}}.\label{BRST}
\eea  

By the above results,  we can write the trace
over
$Ker$ or  $Coker$ as the sum of trace over $\F_{(\a;\b)}^{(l)}$, and 
compute the latter  by using the technique introduced in \cite{Cla73}.
The results are 
\bea
Ch_{Ker_{\F_{(\a;\b)}}}(q;x_1,\cdots,x_N)
&=&\frac{q^{\frac{1}{2}\a(\a-1)}}{\prod_{n=1}^{\infty}(1-q^n)^{N+1}}
\sum_{l=1}^{\infty}(-1)^{l+1}q^{\frac{1}{2}\{l^2+l(2\a-1)\}}\no\\
&&\times \sum_{\{i_1,\cdots,i_N\}\in\Z}
q^{\frac{1}{2}\{i_N^2+i_N(1-2\a-2l)\}}
q^{\frac{1}{2}\D(i_1,\cdots.i_N)}\no\\
&&\times x_1^{2i_1-i_2}x_2^{2i_2-i_1-i_3}\cdots
x_{N-1}^{2i_{N-1}-i_N-i_{N-2}}
x_N^{\a-i_N},\no\\   
Ch_{Coker_{\F_{(\a;\b)}}}(q;x_1,\cdots,x_N)
&=&\frac{q^{\frac{1}{2}\a(\a-1)}}{\prod_{n=1}^{\infty}(1-q^n)^{N+1}}
\sum_{l=1}^{\infty}(-1)^{l+1}q^{\frac{1}{2}\{l^2+l(1-2\a)\}}\no\\
&&\times \sum_{\{i_1,\cdots,i_N\}\in\Z}  
q^{\frac{1}{2}\{i_N^2+i_N(1-2\a+2l)\}}
q^{\frac{1}{2}\D(i_1,\cdots.i_N)}\no\\
&&\times x_1^{2i_1-i_2}x_2^{2i_2-i_1-i_3}\cdots
x_{N-1}^{2i_{N-1}-i_N-i_{N-2}}
x_N^{\a-i_N},\no
\eea
where $\D(i_1,\cdots,i_N)=\sum_{l,l'=1}^{N}\frac{\a_{ll'}\b_{ll'}}
{N-1}\l^l_{i_1,\cdots,i_N}\l^{l'}_{i_1,\cdots,i_N}$  and 
\bea
\left\{
\begin{array}{ll}
\l^l_{i_1,\cdots,i_N}=2i_l-i_{l-1}-i_{l+1},&2\leq l\leq N-1\\
\l^1_{i_1,\cdots,i_N}=2i_1-i_2,&
\l^N_{i_1,\cdots,i_N}=\a-i_N
\end{array}
\right..\label{Weight}
\eea
Similary,  the supercharacters  of $Ker_{\F_{(\a;\b)}}$ and
$Coker_{\F_{(\a;\b)}}$
are   given by   
\begin{enumerate}
\item For $N=2L$:
\bea
Sch_{Ker_{\F_{(\a;\b)}}}(q;x_1,\cdots,x_N)
&=&\frac{(-1)^{\a}
q^{\frac{1}{2}\a(\a-1)}}{\prod_{n=1}^{\infty}(1-q^n)^{N+1}}
\sum_{l=1}^{\infty}(-1)^{l+1}q^{\frac{1}{2}\{l^2+l(2\a-1)\}}\no\\
&&\times  \sum_{\{i_1,\cdots,i_N\}\in\Z} (-1)^{i_N} 
q^{\frac{1}{2}\{i_N^2+i_N(1-2\a-2l)\}}
q^{\frac{1}{2}\D(i_1,\cdots.i_N)}\no\\
&&\times x_1^{2i_1-i_2}x_2^{2i_2-i_1-i_3}\cdots
x_{N-1}^{2i_{N-1}-i_N-i_{N-2}}
x_N^{\a-i_N},\no\\
Sch_{Coker_{\F_{(\a;\b)}}}(q;x_1,\cdots,x_N)
&=&\frac{(-1)^{\a}q^{\frac{1}{2}\a(\a-1)}}{\prod_{n=1}^{\infty}(1-q^n)^{N+1}}
\sum_{l=1}^{\infty}(-1)^{l+1}q^{\frac{1}{2}\{l^2+l(1-2\a)\}}\no\\
&&\times \sum_{\{i_1,\cdots,i_N\}\in\Z} (-1)^{i_N}
q^{\frac{1}{2}\{i_N^2+i_N(1-2\a+2l)\}}
q^{\frac{1}{2}\D(i_1,\cdots.i_N)}\no\\
&&\times x_1^{2i_1-i_2}x_2^{2i_2-i_1-i_3}\cdots
x_{N-1}^{2i_{N-1}-i_N-i_{N-2}}
x_N^{\a-i_N},\no
\eea
\item For $N=2L+1$: 
\bea   
Sch_{Ker_{\F_{(\a;\b)}}}(q;x_1,\cdots,x_N)
&=&-\frac{(-1)^{(L+1)\a}q^{\frac{1}{2}\a(\a-1)}}
{\prod_{n=1}^{\infty}(1-q^n)^{N+1}}
\sum_{l=1}^{\infty}q^{\frac{1}{2}\{l^2+l(2\a-1)\}}\no\\
&&\times \sum_{\{i_1,\cdots,i_N\}\in\Z}(-1)^{i_N}
q^{\frac{1}{2}\{i_N^2+i_N(1-2\a-2l)\}}
q^{\frac{1}{2}\D(i_1,\cdots.i_N)}\no\\
&&\times x_1^{2i_1-i_2}x_2^{2i_2-i_1-i_3}\cdots
x_{N-1}^{2i_{N-1}-i_N-i_{N-2}}
x_N^{\a-i_N},\no\\
Sch_{Coker_{\F_{(\a;\b)}}}(q;x_1,\cdots,x_N)
&=&-\frac{(-1)^{(L+1)\a}q^{\frac{1}{2}\a(\a-1)}}
{\prod_{n=1}^{\infty}(1-q^n)^{N+1}}
\sum_{l=1}^{\infty}q^{\frac{1}{2}\{l^2+l(1-2\a)\}}\no\\
&&\times \sum_{\{i_1,\cdots,i_N\}\in\Z}   (-1)^{i_N}
q^{\frac{1}{2}\{i_N^2+i_N(1-2\a+2l)\}}
q^{\frac{1}{2}\D(i_1,\cdots.i_N)}\no\\
&&\times x_1^{2i_1-i_2}x_2^{2i_2-i_1-i_3}\cdots
x_{N-1}^{2i_{N-1}-i_N-i_{N-2}}
x_N^{\a-i_N}.\no
\eea
\end{enumerate}

Since $\F_{(\a-(N-1);\b+1)}^{(1)}=\F_{(\a;\b)}$ and by  (\ref{BRST}), 
we have
\bea
&&Ch_{Coker_{\F_{(\a-(N-1);\b+1)}}}=Ch_{Ker_{\F_{(\a;\b)}}},~~
Sch_{Coker_{\F_{(\a-(N-1);\b+1)}}}=Sch_{Ker_{\F_{(\a;\b)}}}.\label{C1}
\eea
Relations (\ref{C1}) can also  be checked by using the above explicit
formulae of
the (super)characters. 

\subsection{$\uqsnh$-module structure of $\F_{(\a;\;\b-\frac{1}{N-1}\a)}$}
Set $\l_{\a}=(1-\a)\L_0+\a\L_N$ and 
\bea
&&|\l_{\a}>=|\b,\cdots,\b,\b-\a;-\a>\in \F_{(\a;\b)},~~~\a\in\Z,\no\\
&&|\L_m>=|\b+1,\cdots,\b+1,\b,\cdots,\b;0>\in\F_{(m;\b)},~~m=1,\cdots,N,\no
\eea
The above vectors play the role of the highest weight vectors of 
$\uqsnh$-modules. one can check that 
\bea 
\left\{
\begin{array}{ll}
\eta_0|\l_{\a}>=0,&~~~{\rm for }~\a=0,-1,\cdots \\
\eta_o|\L_m>=0,&~~~{\rm for }~m=1,\cdots,N\\
\eta_0|\l_{\a}>\neq 0,&~~~{\rm for }~\a=1,2,\cdots
\end{array}
\right..\label{HW}
\eea
It follows that the modules 
\bea
&&Coker_{\F_{(\a;\b)}}~~(\a=1,2,\cdots),~~~~
Ker_{\F_{(\a;\b)}}~~(\a=0,-1,-2,\cdots),\no\\
&&Ker_{\F_{(m;\b)}}~~(m=1,2,\cdots,N),\no
\eea
are highest weight $\uqsnh$-modules. Denote them by $\V(\l_{\a})$ and 
$\V(\L_{m})$, respectively. From (\ref{HW}) and (\ref{C1}), 
we have the following identifications of the highest weight 
$\uqsnh$-modules:
\bea 
\V(\l_{\a})&\cong& Ker_{\F_{(\a;\b-\frac{1}{N-1}\a)}}\equiv
Coker_{\F_{(\a-(N-1);\b-\frac{1}{N-1}\a+1)}},~~~{\rm 
for }~\a=0,-1,-2,\cdots,\no\\
&\cong& Coker_{\F_{(\a;\b-\frac{1}{N-1}\a)}}\equiv
Ker_{\F_{(\a+(N-1);\b-\frac{1}{N-1}\a-1)}},~~~{\rm
for }~\a=1,2,\cdots,\label{T31} \\
\V(\L_m)&\cong&Ker_{\F_{(m;\b-\frac{1}{N-1}m)}}\equiv 
Coker_{\F_{(m-(N-1);\b-\frac{1}{N-1}m+1)}},~~~{\rm for}~m=1,\cdots,N.
\label{T32}
\eea
It is easy to see that the vertex operators
(\ref{Vertex-operator2})
also commute (or anti-commute) with $\eta_0$. It follows from 
(\ref{T31})-(\ref{T32}) that  each Fock space
$\F_{(\a;\b-\frac{1}{N-1}\a)}$ is 
decomposed into a direct sum of the highest weight 
$\uqsnh$-modules: 
\bea
 \begin{array}{cccc}
  &Ker && Coker\\
  .&.   && . \\
  .& .&&.\\
F_{(-N;\;\beta+1+\frac{1}{N-1})}=&\V(\lambda_{-N})&\bigoplus&
\V(\l_{-1})\\
     &\phi(z)\uparrow\downarrow \phi^{*}(z)&
     &\phi(z)\uparrow\downarrow \phi^{*}(z)\\
F_{(-N+1;\;\beta+1)}=&\V(\lambda_{-N+1})&\bigoplus& \V(\L_0)\\
     &\phi(z)\uparrow\downarrow \phi^{*}(z)&
     &\phi(z)\uparrow\downarrow \phi^{*}(z)\\
F_{(-N+2;\;\beta+1-\frac{1}{N-1})}=&\V(\l_{-N+2})&\bigoplus& \V(\L_1)\\
     &\phi(z)\uparrow\downarrow \phi^{*}(z)&
     &\phi(z)\uparrow\downarrow \phi^{*}(z)\\
  .&.   && . \\
  .& .&&.\\
F_{(-2;\;\beta+1-\frac{N-3}{N-1})}=&\V(\l_{-2})&\bigoplus& \V(\L_{N-3})\\
     &\phi(z)\uparrow\downarrow \phi^{*}(z)& 
     &\phi(z)\uparrow\downarrow \phi^{*}(z)\\
F_{(-1;\;\beta+1-\frac{N-2}{N-1})}=&\V(\l_{-1})&\bigoplus& \V(\L_{N-2})\\
     &\phi(z)\uparrow\downarrow \phi^{*}(z)& 
     &\phi(z)\uparrow\downarrow \phi^{*}(z)\\
F_{(0;\;\beta)}=&\V(\L_0)&\bigoplus& \V(\L_{N-1})\\
&\phi(z)\uparrow\downarrow \phi^{*}(z)&
     &\phi(z)\uparrow\downarrow \phi^{*}(z)\\
F_{(1;\;\beta-\frac{1}{N-1})}=&\V(\L_1)&\bigoplus& \V(\L_N)\\
     &\phi(z)\uparrow\downarrow \phi^{*}(z)&
     &\phi(z)\uparrow\downarrow \phi^{*}(z)\\
F_{(2;\;\beta-\frac{2}{N-1})}=&\V(\L_2)&\bigoplus& \V(\l_2)\\
     &\phi(z)\uparrow\downarrow \phi^{*}(z)&
     &\phi(z)\uparrow\downarrow \phi^{*}(z)\\
  .&.   && . \\
  .& .&&.\\
F_{(N-2;\;\beta-\frac{N-2}{N-1})}=&\V(\L_{N-2})&\bigoplus&
\V(\l_{N-2})\\
     &\phi(z)\uparrow\downarrow \phi^{*}(z)&
     &\phi(z)\uparrow\downarrow \phi^{*}(z)\\
F_{(N-1;\;\beta-1)}=&\V(\L_{N-1})&\bigoplus& \V(\L_{N-1})\\
     &\phi(z)\uparrow\downarrow \phi^{*}(z)&
     &\phi(z)\uparrow\downarrow \phi^{*}(z)\\
F_{(N;\;\beta-1-\frac{1}{N-1})}=&\V(\L_{N})&\bigoplus& \V(\l_N)\\
     &\phi(z)\uparrow\downarrow \phi^{*}(z)&
     &\phi(z)\uparrow\downarrow \phi^{*}(z)\\
F_{(N+1;\;\beta-1-\frac{2}{N-1})}=&\V(\l_{2})&\bigoplus& \V(\l_{N+1})\\
  .& .&&.\\
 .& .&&.
 \end{array}. \label{D1}
\eea

It is expected that  
$\V(\l_{\a})$ ($\a\in\Z$) and $\V(\L_m)$ ($m=1,2,\cdots,N-1$) are
irreducible highest weight $\uqsnh$-modules  with the highest weights 
$\l_{\a}$ and $\L_m$, respectively. Thus we conjecture that 
\bea 
\V(\l_{\a})=V(\l_{\a}),~~~\V(\L_m)=V(\L_m).\label{IRP}
\eea

\sect{Exchange Relations of Vertex Operators}
In this section, we derive the exchange relations of the type I and type
II bosonized vertex operators of $\uqsnh$. As expected, these 
vertex operators satisfy the graded Faddeev-Zamolodchikov algebra.

\subsection{The R-matrix}
Throughout, we use the abbreviation
\bea
&&(z;x_1,\cdots,x_m)_{\infty}=\prod_{\{n_1,\cdots,n_m\}=0}^{\infty}
(1-zx_1^{n_1}\cdots x_m^{n_m}),\no\\
&&\{z\}_{\infty}\stackrel{def}{=}(z;q^{2(N-1)},q^{2(N-1)})_{\infty}.
\eea

Let $\bar{R}(z) \in End(V\otimes V)$ be the R-matrix of $\uqsnh$, 
\bea
\bar{R}(z)(v_i\otimes v_j)=\sum_{k,l=1}^{2N}\bar{R}^{ij}_{kl}(z)v_k\otimes
v_l~~,
\forall v_i,v_j,v_k,v_l\in V,\label{R}
\eea
\noindent where the matrix elements of $\bar{R}(z)$ are given by 
\begin{eqnarray*}
& &\bar{R}^{i,i}_{i,i}(z)=-1,~~
\bar{R}^{N+1,N+1}_{N+1,N+1}(z)=-\frac{zq^{-1}-q}{zq-q^{-1}},~~
i=1,2,\cdots,N,\\
& &\bar{R}^{ij}_{ij}(z)=\frac{z-1}{zq-q^{-1}},~~~i\neq j,\\
& &\bar{R}^{ji}_{ij}(z)=\frac{q-q^{-1}}{zq-q^{-1}}(-1)^{[i][j]},~~~i<j,\\
& &\bar{R}^{ji}_{ij}(z)=\frac{(q-q^{-1})z}{zq-q^{-1}}(-1)^{[i][j]},
~~~i>j,\\
& &\bar{R}^{ij}_{kl}(z)=0,~~{\rm otherwise}.
\end{eqnarray*} 
Define the R-matrices $R^{(I)}(z)$ and $R^{(II)}(z)$ by
\bea 
R^{(I)}(z)=r(z)\bar{R}(z),~~~~~~~~
R^{(II)}(z)=\bar{r}(z)\bar{R}(z),
\eea
where 
\bea
&&r(z)=z^{\frac{2-N}{N-1}}
\frac{(zq^{2};q^{2(N-1)})_{\infty}(z^{-1}q^{2N-2};q^{2(N-1)})_{\infty}}
{(z^{-1}q^{2};q^{2(N-1)})_{\infty}(zq^{2N-2};q^{2(N-1)})_{\infty}},\no\\
&&\bar{r}(z)=-z^{-\frac{1}{N-1}}
\frac{(zq^{2N-4};q^{2(N-1)})_{\infty}(z^{-1}q^{2N-2};q^{2(N-1)})_{\infty}}
{(z^{-1}q^{2N-4};q^{2(N-1)})_{\infty}(zq^{2N-2};q^{2(N-1)})_{\infty}}.\no
\eea
These R-matrices satisfy the graded Yang-Baxter equation  on 
 $V\otimes V\otimes V$:
\bea
R^{(i)}_{12}(z)R^{(i)}_{13}(zw)R^{(i)}_{23}(w)=
R^{(i)}_{23}(w)R^{(i)}_{13}(zw)R^{(i)}_{12}(z),~~~i=I,II.\no
\eea
Moreover, they  enjoy {\bf (i)} the initial condition 
$R^{(i)}(1)=P,~ i=I,II$, where $P$ is the graded permutation operator;  
{\bf (ii)} the unitarity condition $ 
R^{(i)}_{12}(\frac{z}{w})R^{(i)}_{21}(\frac{w}{z})=1,~~i=I,II$, where 
$R^{(i)}_{21}(z)=PR^{(i)}_{12}(z)P$; {\bf (iii)} 
the crossing-unitarity 
\begin{eqnarray*}
(R^{(i)})^{-1,st_1}(z)\lt((q^{-2\bar\rho}\otimes 1)
R^{(i)}(zq^{2(1-N)})(q^{2\bar\rho}\otimes 1)\rt)
^{st_1}=1,~~~i=I,II,
\end{eqnarray*}
where 
\bea 
q^{2\bar\rho}&\equiv&
diag(q^{2\rho_1},q^{2\rho_2},\cdots,q^{2\rho_{N}},q^{2\rho_{N+1}})\no\\
&=&diag(q^{N-2},q^{N-4},\cdots,q^{-N},q^{-N}).\no
\eea
The various supertranspositions of the R-matrix are given by 
\begin{eqnarray*}
& &(R^{st_1}(z))^{kl}_{ij}=R^{il}_{kj}(z)(-1)^{[i]([i]+[k])},~~~~~
(R^{st_2}(z))^{kl}_{ij}=R^{kj}_{il}(z)(-1)^{[j]([l]+[j])}~,\\
& &(R^{st_{12}}(z))^{kl}_{ij}=R^{ij}_{kl}(z)
(-1)^{([i]+[j])([i]+[j]+[k]+[l])}=R^{ij}_{kl}(z)~~.
\end{eqnarray*}

\subsection{The graded Faddeev-Zamolodchikov algebra}
We now calculate the exchange relations of the type I and type 
II bosonic vertex operators of $\uqsnh$. 
Define 
\begin{eqnarray*}
\oint dzf(z)=Res(f)=f_{-1},~~{\rm for~ a~ formal~ series~ function} 
~ f(z)=\sum_{n\in \Z}f_nz^n.
\end{eqnarray*}
\noindent Then, the Chevalley generators of $\uqsnh$ can be expressed by
the integrals
\begin{eqnarray*}
e_i=\oint dz X^{+,i}(z),~~~~f_i=\oint dz X^{-,i}(z),~i=1,2,\cdots,N.
\end{eqnarray*}
One can also get the integral expressions of the bosonic vertex
operators $\phi(z)$, $\phi^*(z)$, $\psi(z)$ and $\psi^*(z)$ . 
Using these integral expressions and the relations 
given in appendices A and B, we find that   
the bosonic vertex operators defined in (\ref{Vertex-operator2}) satisfy
the graded Faddeev-Zamolodchikov algebra
\bea
&&\phi_j(z_2)\phi_i(z_1)=
\sum_{k,l=1}^{N+1}R^{(I)}(\frac{z_1}{z_2})^{kl}_{ij}
\phi_k(z_1)\phi_l(z_2)(-1)^{[i][j]},\no\\
&&\psi^*_i(z_1)\psi^*_j(z_2)=
\sum_{k,l=1}^{N+1}R^{(II)}(\frac{z_1}{z_2})^{ij}_{kl}
\psi^*_l(z_2)\psi^*_k(z_1)(-1)^{[i][j]},\no\\
&&\psi^*_i(z_1)\phi_j(z_2)=\tau(\frac{z_1}{z_2}) 
\phi_j(z_2)\psi^*_i(z_1)(-1)^{[i][j]},\label{PZF2}
\eea
where 
\bea
\tau(z)=-z^{\frac{2-N}{N-1}}
\frac{(zq;q^{2(N-1)})_{\infty}(z^{-1}q^{2N-3};q^{2(N-1)})_{\infty}}
{(z^{-1}q;q^{2(N-1)})_{\infty}(zq^{2N-3};q^{2(N-1)})_{\infty}}.\no
\eea
By 
\bea 
:e^{-h^*_N(zq^N;\frac{1}{2})+h^*_1(zq;\frac{1}{2})
-h^1(zq^2;\frac{1}{2})-h^2(zq^3;\frac{1}{2})\cdots-
h^N(zq^{N+1};\frac{1}{2})}:=1,\no
\eea
we obtain the first invertibility relations 
\bea
\phi_i(z)\phi^*_j(z)=g^{-1}(-1)^{[i]}\d_{ij},~~~~
\sum_{k=1}^{N+1}(-1)^{[k]}\phi^*_k(z)\phi_k(z)=g^{-1},
\eea
and the second invertibility relations
\bea
\phi^*_i(zq^{2(N-1)})\phi_j(z)=-g^{-1}q^{2\rho_i}\d_{ij},~~~~
\sum_{k=1}^{N+1}q^{-2\rho_k}\phi_k(z)\phi^*_k(zq^{2(N-1)})=-g^{-1},
\label{IV}
\eea
where $g=e^{\frac{\sqrt{-1}\pi N}{2(N-1)}}\frac{(q^2;q^{2(N-1)})_{\infty}}
{(q^{2(N-1)};q^{2(N-1)})_{\infty}}$. 
Using the fact that $\eta_0\xi_0$ is a projection operator, 
we can make the following identifications:
\bea
&&\Phi_i(z)=\eta_0\xi_0\phi_i(z)\eta_0\xi_0,~~~~~
\Phi^*_i(z)=\eta_0\xi_0\phi^*_i(z)\eta_0\xi_0,\no\\
&&\Psi_i(z)=\eta_0\xi_0\psi_i(z)\eta_0\xi_0,~~~~~
\Psi^*_i(z)=\eta_0\xi_0\psi^*_i(z)\eta_0\xi_0.\label{V2}
\eea
Set 
\bea 
\mu_{\a}=
\left\{
\begin{array}{ll}
\L_{\a},&\a=0,1,\cdots,N\\
\l_{\a-(N-1)},&{\rm for}~\a>N\\
\l_{\a},&{\rm for}~\a<0
\end{array}
\right..\label{M1}
\eea
It is easy to see that  the vertex operators $\phi(z)$, $\phi^*(z)$,
$\psi(z)$ and 
$\psi^*(z)$ commute (or anti-commute) with the BRST charge $\eta_0$. 
It follows from   (\ref{D1}) and (\ref{IRP}) that   
the vertex operators (\ref{V2}) intertwine all the
level-one irreducible highest weight $\uqsnh$-modules $V(\mu_{\a})$
($\a\in\Z$) as follows
\bea
& &\Phi(z) :
 V(\mu_{\a}) \longrightarrow V(\mu_{\a-1})
\otimes V_{z} ,\ \
\Phi^{*}(z) :
 V(\mu_{\a}) \longrightarrow V(\mu_{\a+1})\otimes
V_{z}^{*S},\no\\
& &\Psi(z) :
 V(\mu_{\a}) \longrightarrow V_{z}\otimes V(\mu_{\a-1})
,\ \
\Psi^{*}(z) :
 V(\mu_{\a}) \longrightarrow V_{z}^{*S}\otimes
V(\mu_{\a+1}).\label{Mu1}
\eea
>From (\ref{PZF2}), we have 
\bea
&&\Phi_j(z_2)\Phi_i(z_1)=
\sum_{k,l=1}^{N+1}R^{(I)}(\frac{z_1}{z_2})^{kl}_{ij}
\Phi_k(z_1)\Phi_l(z_2)(-1)^{[i][j]},\no\\
&&\Psi^*_i(z_1)\Psi^*_j(z_2)=
\sum_{k,l=1}^{N+1}R^{(II)}(\frac{z_1}{z_2})^{ij}_{kl}
\Psi^*_l(z_2)\Psi^*_k(z_1)(-1)^{[i][j]},\no\\
&&\Psi^*_i(z_1)\Phi_j(z_2)=\tau(\frac{z_1}{z_2})
\Phi_j(z_2)\Psi^*_i(z_1)(-1)^{[i][j]}.\label{ZF2}
\eea
Moreover, we have  the following 
invertibility relations:
\bea
&&\Phi_i(z)\Phi^*_j(z)=g^{-1}(-1)^{[i]}\d_{ij}id_{V(\mu_{\a})},\no\\
&&\sum_{k=1}^{N+1}(-1)^{[k]}\Phi^*_k(z)\Phi_k(z)=
g^{-1}id_{V(\mu_{\a})},\no\\
&&\Phi^*_i(zq^{2(N-1)})\Phi_j(z)=-g^{-1}q^{2\rho_i}\d_{ij}
id_{V(\mu_{\a})},\no\\
&&\sum_{k=1}^{N+1}q^{-2\rho_k}\Phi_k(z)\Phi^*_k(zq^{2(N-1)})=-g^{-1}
id_{V(\mu_{\a})}.\label{Inv}
\eea

\sect{Multi-component super  $t$-$J$ model }
In this section, we give a mathematical definition of the 
multi-component  super $t$-$J$ model on an infinite lattice.

\subsection{Space of states}
By means of the R-matrix (\ref{R}) of $\uqsnh$,
one  defines  a spin chain  model, referred to as the multi-component 
super $t-J$ model,  on the infinte lattice $\cdots
\otimes V\otimes V\otimes V\cdots$. Let $h$ be the operator on
$V\otimes V$ such that
\begin{eqnarray*}
& &P\bar{R}(\frac{z_1}{z_2})=1+u h+\cdots,\ \ \ \ \ \ \ \ u
\longrightarrow 0 ,\\
& &\ \ \ \ \  P:{\rm the\ \  graded\ \ permutation \ \ operator },\ \
 e^{u}\equiv \frac{z_1}{z_2}.
\end{eqnarray*}
The Hamiltonian $H$ of this model is given  by
\begin{eqnarray}
H=\sum_{l\in Z}h_{l+1,l}.
\end{eqnarray}
$H$ acts formally on the infinite tensor product,
\beq
\cdots V\otimes V\otimes V\cdots.\label{vvv}
\eeq
It can be easily checked that
\begin{eqnarray*}
[U'_q(\widehat{sl}(N|1)), H]=0,
\end{eqnarray*}
where $U'_q[\widehat{sl}(N|1)]$ is the subalgebra of $\uqsnh$ with the
derivation operator $d$ being dropped. So $U'_q[\widehat{sl}(N|1)]$ plays
the role of infinite dimensional {\it non-abelian symmetry} of the 
multi-component  super $t$-$J$ model on the infinite lattice.

>From the intertwining relation (\ref{Mu1}), one have the following 
composition of the type I vertex operators:
\bea 
V(\mu_{\a})\stackrel{\Phi(1)}{\longrightarrow}V(\mu_{\a-1})\otimes V
\stackrel{\Phi(1)\otimes id}{\longrightarrow}V(\mu_{\a-1})\otimes V
\otimes V\stackrel{\Phi(1)\otimes
id\otimes id }{\longrightarrow}\cdots\longrightarrow W_{l},\label{CM}
\eea
where $W_l\stackrel{def}{=}\cdots\otimes V\otimes V$, i.e the
left half-infinite tensor product. We conjecture that such a composition 
converges to a map :
\bea 
i: V(\mu_{\a})\longrightarrow W_l.\no
\eea
Such a map $i$ satisfies 
$i(xv)=\D^{(\infty)}(x)i(v)$, $ x\in \uqsnh$ and 
$v \in V(\mu_{\a})$.
Following \cite{Dav93}, we could replace the infinite tensor product
(\ref{vvv}) by the level-zero $\uqsnh$-module,
\begin{eqnarray*}
{\it F}_{\a\a'}={\rm Hom}(V(\mu_\a),V(\mu_{\a'}))\cong V(\mu_\a)\otimes
  V(\mu_{\a'})^{*},
\end{eqnarray*}
where $V(\mu_\a)$ is level-one irreducible highest weight
$\uqsnh$-module and $V(\mu_{\a'})^{*}$ is the dual module of
$V(\mu_{\a'})$.
By (\ref{D1}), this homomorphism can be realized by applying the type I
vertex
operators repeatedly.
So we shall make the (hypothetical) identification:
\begin{eqnarray*}
``{\rm \bf the \ \ space \ \ of \ \ physical \ \ states } " =
\bigoplus_{\a,\a'\in \Z} V (\mu_\a)\otimes V(\mu_{\a'})^{*}.
\end{eqnarray*}
Namely, we take
\begin{eqnarray*}
F\equiv End(\bigoplus_{\a\in \Z}V(\mu_\a))\cong
 \bigoplus_{\a,\a'\in{\bf Z}} F_{\a\a'}
\end{eqnarray*}
as the space of states of the multi-component super $t$-$J$
model on the infinite lattice.
The left action of $\uqsnh$ on ${\it F}$ is defined by
\begin{eqnarray*}
x.f=\sum x_{(1)}\circ f\circ S(x_{(2)})(-1)^{[f][x_{(2)}]},~~~
 \forall x\in\uqsnh,~f\in F,
\end{eqnarray*}
where we have used notation $\Delta(x)=\sum x_{(1)}
\otimes x_{(2)}$.
Note that
$F_{\a\a}$ has the unique canonical element $id_{V(\mu_\a)}$. We call
it the vacuum \cite{Jim94} and denote it by
$|vac>_{\a}$.
\subsection{Local structure and local operators}
Following Jimbo et al \cite{Jim94}, we use the type I
vertex operators and their variants to incorporate the local structure
into the
space of physical states ${\it F}$, that is to formulate the action
of local operators of the multi-component super $t$-$J$
model on the infinite tensor product (\ref{vvv}) in terms of their actions
on $F_{\a\a'}$.   

Using the isomorphisms 
\bea
\Phi(1)&:&~V(\mu_\a)\longrightarrow V(\mu_{\a-1})\otimes V,\no\\
\Phi^{*,st}(q^{2(N-1)})&:&~V\otimes V(\mu_\a)^*\longrightarrow
  V(\mu_{\a-1})^*,
\eea
were $st$ is the supertransposition on the quantum space,
we have the following identification:
\begin{eqnarray*}
V(\mu_\a)\otimes V(\mu_{\a'})^{*}\stackrel{\sim}{\rightarrow}
V(\mu_{\a-1})\otimes V\otimes V(\mu_{\a'})^{*}\stackrel{\sim}{\rightarrow}
V(\mu_{\a-1})\otimes V(\mu_{\a'-1})^{*}.
\end{eqnarray*}
The resulting isomorphism can be identified with
the super translation (or shift) operator defined by
\begin{eqnarray*}
T=-g\sum_i\Phi_i(1)\otimes
\Phi_i^{*,st}(q^{2(N-1)})(-1)^{[i]}q^{-2\rho_i}.
\end{eqnarray*}
Its inverse is given by
\begin{eqnarray*}
T^{-1}=g\sum_i\Phi_i^{*}(1)\otimes \Phi_i^{st}(1).
\end{eqnarray*}

Thus we can define the local operators on $V$ as operators on
$F_{\a\a'}$ \cite{Jim94}. Let us label the tensor components from the
middle
as
$1,2,\cdots$ for the left half and as $0,-1,-2,\cdots$ for
the right half. The operators acting on the site 1  are defined by
\begin{eqnarray}
E_{ij}\stackrel{def}{=}E^{(1)}_{ij}=
g\Phi^{*}_i(1)\Phi_j(1)(-1)^{[j]}\otimes id.
\end{eqnarray}
More generally we set
\begin{eqnarray}
\label{L}
E^{(n)}_{ij}=T^{-(n-1)}E_{ij}T^{n-1}\ \ \ \ (n\in Z).
\end{eqnarray} 
Then, from the invertibility relations of the type I vertex operators
of $\uqsnh$, we can show that 
the local operators $E^{(n)}_{ij}$ acting on
${\it F}_{\a\a'}$  satisfy the following relations:
\begin{eqnarray*}
E^{(m)}_{ij}E^{(n)}_{kl}=\left\{
  \begin{array}{ll}
  \delta_{jk}E^{(n)}_{il} & {\rm if}\ \ \ m=n\\
  (-1)^{([i]+[j])([k]+[l])}E^{(n)}_{kl}E^{(m)}_{ij}&{\rm if}\ \ \ m\ne n
  \end{array}   
 \right. .
\end{eqnarray*}

This  result implies that the local
operators $E^{(n)}_{ij}$ are nothing but the $U_q[sl(N|1)]$ 
 generators  acting on the n-th component of 
$\cdots \otimes V\otimes V\otimes \cdots$. They include 
all the local operators in the multi-component  super $t-J$
model \cite{Jim94}.  

As is expected from the physical point of view, the vacuum vectors
$|vac>_{\a}$ are super-translationally invariant and singlets (i.e. 
they belong
to the trivial representation of $\uqsnh$):
\begin{eqnarray*}
T|vac>_\a=|vac>_{\a-1},~~~~~
x.|vac>_\a=\epsilon(x)|vac>_\a,~~\forall x\in \uqsnh.
\end{eqnarray*}
This is proved as follow. 
 Let $u^{(\a)}_l$ ($u^{*(\a)}_l$) be a basis
vectors
of
$V(\mu_\a)\ \  (V(\mu_\a)^{*})$  and
\begin{eqnarray*}
|vac>_\a\stackrel{def}{=}id_{V(\mu_\a)}=\sum_l u^{(\a)}_l\otimes
u^{*(\a)}_l.
\end{eqnarray*}
Then
\bea
T|vac>_\a=-g\sum_{m,l}q^{-2\rho_m}\Phi_m(1)u^{(\a)}_l\otimes \Phi^{*,st}_m  
(q^{2(N-1)})u^{*(\a)}_l(-1)^{[m]+[l][m]}.\no
\eea
We want to show $T|vac>_\a=|vac>_{\a-1}$. This
is equivalent to proving
\begin{eqnarray*}
-g\sum_{m,l}q^{-2\rho_m}\Phi_m(1)u^{(\a)}_l\,\Phi^{*,st}_m(q^{2(N-1)})\cdot
  u^{*(\a)}_l(v)(-1)
  ^{[m]+[l][m]}=v,~~~\forall v\in V(\mu_{\a-1}).
\end{eqnarray*}
Now
\bea
 l.h.s&=&-g\sum_{m,l}q^{-2\rho_m}\Phi_m(1)u^{(\a)}_l u^{*(\a)}_l\lt
(\Phi^{*}_m(q^{2(N-1)})^{st})^{st}v\rt)(-1)^{[m]}\no\\
&=&-g\sum_{m,l}q^{-2\rho_m}\Phi_m(1)u^{(\a)}_l u^{*(\a)}_l
  (\Phi^{*}_m(q^{2(N-1)})v)\no\\
&=&-g\sum_{m}q^{-2\rho_m}\Phi_m(1)\Phi^{*}_m(q^{2(N-1)})v=v ,\no
\eea
where we have used $(\Phi^{*}_m(z)^{st})^{st}
=\Phi^{*}_m(z)(-1)^{[m]}$ and (\ref{Inv}). As to the second equation, we
have   
\bea
x\cdot|vac>_\a&=&\sum x_{(1)}u^{(\a)}_l\otimes x_{(2)}u^{*(\a)}_l
  (-1)^{[l][x_{(2)}]}\no\\
&=&\sum x_{(1)}u^{(\a)}_l\otimes
\pi_{V(\mu_\a)^*}(x_{(2)})_{ml}u^{*(\a)}_m
  (-1)^{[l][x_{(2)}]}\no\\
&=&\sum x_{(1)}u^{(\a)}_l\otimes \pi_{V(\mu_\a)}
  (S(x_{(2)}))_{lm}u^{*(\a)}_m\no\\
&=&\sum x_{(1)}\pi_{V(\mu_\a)}(S(x_{(2)}))_{lm}u^{(\a)}_l\otimes
u^{*(\a)}_m
  \no\\
&=&\sum x_{(1)}S(x_{(2)})u^{(\a)}_m\otimes u^{*(\a)}_m 
  =\e(x)|vac>_\a.\no
\eea
This completes the proof.

For any local operator $O\in F$, its vacuum expectation value
is defined  by
\bea
{}_\a<vac|O|vac>_\a\stackrel{def}{=}
\frac{tr_{V(\mu_\a)}(q^{-2\rho}O)}
  {tr_{V(\mu_\a)}(q^{-2\rho})}
=\frac{tr_{V(\mu_\a)}(q^{-2(N-1)d-2h_{\bar{\rho}}}O)}
  {tr_{V(\mu_\a)}(q^{-2(N-1)d-2h_{\bar{\rho}}})},
\eea
where  
\bea
2h_{\bar{\rho}}=\sum_{l=1}^{N}l(N-1-l)h_l.\no
\eea
We shall
denote the correlator ${}_\a<vac|O|vac>_\a$ by $<O>_\a$.

\sect{Correlation functions}

The aim of  this section is to calculate $<E_{mn}>_{\a}$.
The generalization to the calculation of the  multi-point functions
is straightforward.

Set
\begin{eqnarray*}
P^m_n(z_1,z_2|q|\a)=
\frac{tr_{V(\mu_\a)}(q^{-2(N-1)d-2h_{\bar{\rho}}}
\Phi^{*}_m(z_1)\Phi_n(z_2))}
{tr_{V(\mu_\a)}(q^{-2(N-1)d-2h_{\bar{\rho}}})} ,
\end{eqnarray*}
then $<E_{mn}>_{\a}=P^m_n(z,z|q|\a)$. By (\ref{M1}), it is sufficient to
calculate
\begin{eqnarray}
F^{(\alpha)}_{mn}(z_1,z_2)=
\frac{tr_{F_{(\alpha;\beta-\alpha)}}(q^{-2(N-1)d-2h_{\bar{\rho}}}
\phi^{*}_m(z_1)\phi_n(z_2)\eta_0\xi_0)} 
{tr_{F_{(\alpha;\beta-\alpha)}}(q^{-2(N-1)d-2h_{\bar{\rho}}}\eta_0\xi_0)}.
\end{eqnarray}
Using the Clavelli-Shapiro technique \cite{Cla73}, we
get
\begin{eqnarray*}   
F^{(\alpha)}_{mn}(z_1,z_2)=\frac{\delta_{mn}}{\chi_{\a}}
F^{(\alpha)}_m(z_1,z_2)\equiv 
\frac{\delta_{mn}}{\chi_{\a}}
\sum_{l=1}^\infty(-1)^{l+1}F^{(\a)}_{m,-l}(z_1,z_2),
\end{eqnarray*}
where
\bea 
\chi_{\a}&=&Ch_{Ker_{\F_{(\a;\b)}}}(q^{2(N-1)};q^{-(N-2)},\cdots,
q^{-l(N-1-l)},\cdots,q^{N}),\no\\
F^{(\a)}_{m,l}(z_1,z_2)&=&-e^{\frac{\sqrt{-1}\pi N}{2(N-1)}}C^*_1
C^*_N (C_1)^{N-1}(C_{N+1})^2\no\\
&&(z_1q)^{\frac{1}{N-1}}
\frac{\{\frac{z_1}{z_2}q^{2(N-1)}\}_{\infty}\{\frac{z_2}{z_1}
q^{2(N-1)}\}_{\infty}}
{\{\frac{z_1}{z_2}q^{2N}\}_{\infty}\{\frac{z_2}{z_1}q^{2N}\}_{\infty}}
 \oint dw_1\cdots\oint dw_N \no\\
&&\times 
\{\prod_{k=1}^{m-1}\frac{(1-q^2)}
{qw_{k-1}(\frac{w_{k}}{w_{k-1}}q;q^{2(N-1)})_{\infty}
(\frac{w_{k-1}}{w_k}q;q^{2(N-1)})_{\infty}}\}
\no\\
&&\times \frac{1}{w_{m-1}(\frac{w_{m}}{w_{m-1}}q;q^{2(N-1)})_{\infty}
(\frac{w_{m-1}}{w_m}q^{2N-1};q^{2(N-1)})_{\infty}}\no\\
&&\times\{\prod_{k=m+1}^{N}\frac{(1-q^2)}
{w_k(\frac{w_{k}}{w_{k-1}}q;q^{2(N-1)})_{\infty}
(\frac{w_{k-1}}{w_k}q;q^{2(N-1)})_{\infty}}\}
\no\\
&&\times\sum_{\{i_1,\cdots,i_N\}\in\Z}I^{(a,\;l)}_{i_1,\cdots,i_N}
(z_1,z_2|w_1,\cdots,w_N)\no\\
&&\times \lt\{
\frac{(\frac{z_2}{w_N}q^{N-1})^{l-\a+i_N}}
{w_Nq(\frac{z_2}{w_N}q^{N-1};q^{2(N-1)})_{\infty}
(\frac{w_N}{z_2}q^{N-1};q^{2(N-1)})_{\infty}}
\rt.\no\\
&&~~\lt. +
\frac{(\frac{z_2}{w_N}q^{N+1})^{l-\a+i_N}}
{z_2q^{N}(\frac{z_2}{w_N}q^{3N-1};q^{2(N-1)})_{\infty}
(\frac{w_N}{z_2}q^{-N-1};q^{2(N-1)})_{\infty}}\rt\},\no\\
&&{\rm for}~~m=1,\cdots,N,\no\\
\no\\
F^{(\a)}_{N+1,l}(z_1,z_2)&=&e^{\frac{\sqrt{-1}\pi N}{2(N-1)}}C^*_1
C^*_N (C_1)^N(C_{N+1})^2
(z_1q)^{\frac{1}{N-1}}\frac{\{\frac{z_1}{z_2}q^{2(N-1)}\}_{\infty}
\{\frac{z_2}{z_1}  
q^{2(N-1)}\}_{\infty}}
{\{\frac{z_1}{z_2}q^{2N}\}_{\infty}\{\frac{z_2}{z_1}q^{2N}\}_{\infty}
}\no\\
&&\times \oint dw_1\cdots\oint dw_N
\{\prod_{k=1}^{N}\frac{(1-q^2)}
{qw_{k-1}(\frac{w_{k}}{w_{k-1}}q;q^{2(N-1)})_{\infty}
(\frac{w_{k-1}}{w_k}q;q^{2(N-1)})_{\infty}}\}
\no\\
&&\times \frac{1}{w_{N}(\frac{z_2}{w_N}q^{N+1};q^{2(N-1)})_{\infty}
(\frac{w_{N}}{z_2}q^{N-1};q^{2(N-1)})_{\infty}}\no\\
&&\times\sum_{\{i_1,\cdots,i_N\}\in\Z}
I^{(a,\;l)}_{i_1,\cdots,i_N}(z_1,z_2|w_1,\cdots,w_N)\no\\
&&~~\times \partial_{w_N}\{
\frac{(\frac{z_2}{w_N}q^{N})^{l-\a+i_N}}
{w_N(\frac{z_2}{w_N}q^{N};q^{2(N-1)})_{\infty}
(\frac{w_N}{z_2}q^{N-2};q^{2(N-1)})_{\infty}}\}.\no
\eea
In the above equations, $w_0\equiv z_1q$, and
\bea 
I^{(a,\;l)}_{i_1,\cdots,i_N}(z_1,z_2|w_1,\cdots,w_N)&=&
q^{(N-1)\a(\a-1)}(z_1q)^{i_1-\frac{\a}{N-1}}(z_2q^N)^{\frac{N}{N-1}\a-i_N}
\no\\
&&\times 
q^{(N-1)\{l^2+l(1-2\a)+i_N^2+i_N(1-2\a+2l)+\D(i_1,\cdots.i_N)\}}\no\\
&&\times \prod_{k=1}^{N}(w_kq^{k(N-1-k)})^{-\l^k_{i_1,\cdots,i_N}},\no
\eea
\bea
&&C^*_1=\frac{\{q^{2N}\}_{\infty}}{\{q^{4N-4}\}_{\infty}},~~~~~~~~~
C^*_N=\frac{\{q^{4N-2}\}_{\infty}}{\{q^{2(N-1)}\}_{\infty}},\no\\
&&C_1=(q^{2(N-1)};q^{2(N-1)})_{\infty}(q^{2N};q^{2(N-1)})_{\infty},~~~
C_{N+1}=(q^{2(N-1)};q^{2(N-1)})_{\infty}.\no
\eea

We now derive the difference equations satisfied by these one-point
functions. 
Noticing that
\begin{eqnarray*}
& &x^{d}\phi_i(z)x^{-d}=\phi_i(zx^{-1}),\ \
x^{d}\phi^{*}_i(z)x^{-d}=\phi^{*}_i(zx^{-1}),\\
& &x^{d}\psi_i(z)x^{-d}=\psi_i(zx^{-1}),\ \
x^{d}\psi^{*}_i(z)x^{-d}=\psi^{*}_i(zx^{-1}),\\
& &x^{d}\eta_0x^{-d}=\eta_0,\ \
x^{d}\xi_0x^{-d}=\xi_0,
\end{eqnarray*}
we get the difference equations
\begin{eqnarray*}
F^{(\alpha)}_m(z_1,z_2q^{2(N-1)})=q^{-2\rho_m}\,\sum_{k}R(z_2,z_1)^{km}_{mk}
F^{(\alpha-1)}_k(z_1,z_2)
 \,(-1)^{[m]+[k]+[m][k]}.
\end{eqnarray*}
Since $\a\in{\bf Z}$, it is easily seen that this is a set of infinite
number of difference equations.

\vskip.3in
\subsection*{ Acknowledgements}

This work has been  financially supported by the Australian Research 
Council  Large, Small and QEII Fellowship grants. 
W.-L. Yang would like to thank Y.-Z. Zhang, and 
the department of Mathematics, the University of Queensland,
for their kind  hospitality. W.-L. Yang has also been partially supported
by the National Natural Science Foundation of China.

\subsection*{Appendix A}
In this appendix, we give the normal ordered relations of the fundamental 
bosonic fields:

\begin{eqnarray*}
&&:e^{h^i(z;\b_1)}::e^{h^j(w;\b_2)}:=z^{a_{ij}}(1-\frac{w}{z}
q^{\b_1+\b_2})^{a_{ij}}:e^{h^i(z;\b_1)+h^j(w;\b_2)}:,~~~i\neq j,\\
&&:e^{h^i(z;\b_1)}::e^{h^i(w;\b_2)}:=z^2(1-\frac{w}{z}q^{\b_1+\b_2-1})
(1-\frac{w}{z}q^{\b_1+\b_2+1}) 
:e^{h^i(z;\b_1)+h^i(w;\b_2)}:,~~~i\neq N,\\
&&:e^{h^N(z;\b_1)}::e^{h^N(w;\b_2)}:=:e^{h^N(z;\b_1)+h^N(w;\b_2)}:,\\
&&:e^{h^i(z;\b_1)}::e^{h^*_j(w;\b_2)}:=z^{\d_{ij}}(1-\frac{w}{z}
q^{\b_1+\b_2})^{\delta_{ij}}:e^{h^i(z;\b_1)+h^*_j(w;\b_2)}:,\\      
&&:e^{h^*_i(z;\b_1)}::e^{h^{j}(w;\b_2)}:=z^{\d_{ij}}
(1-\frac{w}{z}q^{\b_1+\b_2})^{\delta_{ij}}
:e^{h^*_i(z;\b_1)+h^{j}(w;\b_2)}:,\\
&&:e^{h^*_N(z;\b_1)}::e^{h^*_N(w;\b_2)}:=
z^{-\frac{N}{N-1}}\frac{(\frac{w}{z}q^{\b_1+\b_2+2N-1};q^{2(N-1)})}
{(\frac{w}{z}q^{\b_1+\b_2-1};q^{2(N-1)})}
:e^{h^*_N(z;\b_1)+h^*_N(w;\b_2)}:,\\
&&:e^{h^*_1(z;\b_1)}::e^{h^*_1(w;\b_2)}:=
z^{\frac{N-2}{N-1}}\frac{(\frac{w}{z}q^{\b_1+\b_2+1};q^{2(N-1)})}
{(\frac{w}{z}q^{\b_1+\b_2+2N-3};q^{2(N-1)})}
:e^{h^*_1(z;\b_1)+h^*_1(w;\b_2)}:,\\
&&:e^{h^*_1(z;\b_1)}::e^{h^*_N(w;\b_2)}:=
z^{-\frac{1}{N-1}}\frac{(\frac{w}{z}q^{\b_1+\b_2+N};q^{2(N-1)})}
{(\frac{w}{z}q^{\b_1+\b_2+N-2};q^{2(N-1)})}
:e^{h^*_1(z;\b_1)+h^*_N(w;\b_2)}:,\\
&&:e^{h^*_N(z;\b_1)}::e^{h^*_1(w;\b_2)}:=
z^{-\frac{1}{N-1}}\frac{(\frac{w}{z}q^{\b_1+\b_2+N};q^{2(N-1)})}
{(\frac{w}{z}q^{\b_1+\b_2+N-2};q^{2(N-1)})}
:e^{h^*_N(z;\b_1)+h^*_1(w;\b_2)}:,\\
&&:e^{c(z;\b_1)}::e^{c(w;\b_2)}:=z(1-\frac{w}{z}q^{\b_1+\b_2})
:e^{c(z;\b_1)+c(w;\b_2)}:,\\
\end{eqnarray*}
\noindent where $a_{ij}$ is the Cartan-matrix of
$\widehat{sl(N|1)}$  and $i,j=1,2,\cdots,N$.

\subsection*{Appendix B}
By means of the bosonic realization (\ref{A2}) of $\uqsnh$, 
the integral expressions of the bosonized vertex operators 
(\ref{Vertex-operator2}) and 
the technique given in \cite{AJMP}, one can check the following
relations

\begin{itemize}

\item For the type I vertex operators:

\begin{eqnarray*}
&&[\phi_k(z),f_l]=0 ~{\rm if}~k\neq l,l+1,~~
~[\phi_{l+1}(z),f_l]_{q^{\nu_{l+1}}}=
\nu_l\phi_l(z)(-1)^{[f_l]([v_l]+[v_{l+1}])},\\
&&[\phi_l(z),f_l]_{q^{-\nu_l}}=0,~~  
[\phi_{l}(z),e_l]=q^{h_l}\phi_{l+1}(z)(-1)^{[e_l]([v_l]+[v_{l+1}])},\\
&&[\phi_k(z),e_l]=0 ~{\rm if}~k\neq l,~~
~~q^{h_l}\phi_l(z)q^{-h_l}=q^{-\nu_l}\phi_l(z),\\   
&&q^{h_l}\phi_k(z)q^{-h_l}=\phi_k(z)~{\rm if}~k\neq l,l+1,~~
q^{h_l}\phi_{l+1}(z)q^{-h_l}=q^{\nu_{l+1}}\phi_{l+1}(z),\\
\\
&&[\phi^*_k(z),f_l]=0 ~{\rm if}~k\neq l,l+1,~~
[\phi^*_{l+1}(z),f_l]_{q^{-\nu_{l+1}}}=0,~~\\
&&[\phi^*_k(z),e_l]=0 ~{\rm if}~k\neq l+1,~~
[\phi^*_{l+1}(z),e_l]=-\nu_l\nu_{l+1}q^{h_l-\nu_l}\phi^*_{l}(z)
(-1)^{[e_l]([v_l]+[v_{l+1}])},\\
&&[\phi^*_{l}(z),f_l]_{q^{\nu_l}}=-\nu_lq^{\nu_l}\phi^*_{l+1}(z)
(-1)^{[f_l]([v_l]+[v_{l+1}])},~~
q^{h_l}\phi^*_l(z)q^{-h_l}=q^{\nu_l}\phi^*_l(z),\\
&&q^{h_l}\phi^*_k(z)q^{-h_l}=\phi^*_k(z)~{\rm if}~k\neq l,l+1,~~
q^{h_l}\phi^*_{l+1}(z)q^{-h_l}=q^{-\nu_{l+1}}\phi^*_{l+1}(z).
\\
\end{eqnarray*}
\item For the type II vertex operators:
\begin{eqnarray*}
&&[\psi_k(z),e_l]=0 ~{\rm if}~k\neq l,l+1,~~
[\psi_{l+1}(z),e_l]_{q^{-\nu_{l+1}}}=0,~~
[\psi_{l}(z),e_l]_{q^{\nu_l}}=\psi_{l+1}(z),\\
&&[\psi_k(z),f_l]=0 ~{\rm if}~k\neq l+1,~~
[\psi_{l+1}(z),f_l]=\nu_lq^{-h_l}\psi_{l}(z),~~\\
&&q^{h_l}\psi_l(z)q^{-h_l}=q^{-\nu_l}\psi_l(z),~~
q^{h_l}\psi_{l+1}(z)q^{-h_l}=q^{\nu_{l+1}}\psi_{l+1}(z),\\
&&q^{h_l}\psi_k(z)q^{-h_l}=\psi_k(z)~{\rm if}~k\neq l,l+1,~~\\
\\
&&[\psi^*_k(z),e_l]=0 ~{\rm if}~k\neq l,l+1,~~
[\psi^*_{l}(z),e_l]_{q^{-\nu_l}}=0,~~\\
&&[\psi^*_k(z),f_l]=0 ~{\rm if}~k\neq l,~~   
[\psi^*_{l}(z),f_l]=-\nu_lq^{-h_l+\nu_l}\psi^*_{l+1}(z),~~\\
&&[\psi^*_{l+1}(z),e_l]_{q^{\nu_{l+1}}}=-\nu_l\nu_{l+1}q^{-\nu_l}
\psi^*_{l}(z),~~
q^{h_l}\psi^*_l(z)q^{-h_l}=q^{\nu_l}\psi^*_l(z),\\
&&q^{h_l}\psi^*_k(z)q^{-h_l}=\psi^*_k(z)~{\rm if}~k\neq l,l+1,~~
q^{h_l}\psi^*_{l+1}(z)q^{-h_l}=q^{-\nu_{l+1}}\psi^*_{l+1}(z).
\end{eqnarray*}
\end{itemize}


\end{document}